\documentclass{article}
\usepackage{amsmath,amssymb}
\newtheorem{df}{Definition}[section]
\newtheorem{thm}[df]{Theorem}
\newtheorem{lem}[df]{Lemma}
\newtheorem{prop}[df]{Proposition}

\title{Multiresolution on $n$--dimensional spheres}
\author{I. Iglewska-Nowak\footnote{West Pomeranian University of Technology in Szczecin, School of Mathematics, al. Piast\'ow 17, 70--310 Szczecin, Poland}}

\begin{document}

\maketitle

\bibliographystyle{amsplain}

\begin{abstract}In the present paper, multiscale systems of polynomial wavelets on an $n$--dimensional sphere are constructed. Scaling functions and wavelets are investigated, their reproducing and localization properties and positive definiteness are examined. Decomposition and reconstruction algorithms for the wavelet transform are presented. Formulae for variances in space and momentum domain, as well as for the uncertainty product, of zonal functions over $n$--dimensional spheres are derived and applied to the scaling functions.\end{abstract}

\begin{bfseries}Keywords:\end{bfseries} $n$--dimensional spheres, polynomial wavelets, multiresolution analysis, uncertainty product \\
\begin{bfseries}MSC2010:\end{bfseries} 42C40, 42C05

\section{Introduction}

In the past two decades, several wavelet constructions over $n$--dimensio\-nal sphe\-res were developed, cf., e.g., \cite{AVn,EBCK09,sE11}, all of them being a generalization of the two--dimensio\-nal case. Our aim in this paper is to propose another approach, based on the book chapter of Conrad and Prestin~\cite{MRes}. Wavelets here are polynomials in spherical harmonics, a feature that on the one hand results in an oscillatory behavior of the wavelets themselves, on the other hand, it allows us to define a multiresolution analysis of sampling spaces over equiangular grids. Other MRA constructions for the $2$--sphere can be found, e.g., in~\cite{FGS-book,BP02,FS06,SMAN06}. To the author's best knowledge, the present research is the first attempt to define an MRA over~$\mathcal S^n$. Although the generalization of~\cite{MRes} is straightforward, it demanded much carefulness, especially in the generalization of the Clenshaw--Curtis quadrature in Lemma~\ref{lem:quadrature} and in the the computation of constants in Theorem~\ref{thm:wsf_properties} and Theorem~\ref{thm:ww_proprties}. (and, consequently, constants for reconstruction and decomposition algorithms of Subsection~\ref{subs:algorithms}), since the addition theorem has been reprinted with a false constant since~\cite{LV} and it has been corrected only recently~\cite[formula~(2)]{IIN14CWT}. (Note that the mistake was caused by a change of normalization convention for hyperspherical harmonics, from
$$
\int_{\mathcal S^n}\left|Y_l^k(x)\right|^2d\sigma(x)=1
$$
to
$$
\frac{1}{\Sigma_n}\int_{\mathcal S^n}\left|Y_l^k(x)\right|^2d\sigma(x)=1.
$$
For the two--dimensional sphere, one still uses the first one such that the formulae presented in the numerous papers about analysis over~$\mathcal S^2$ remain unaltered.) We also derive a formula for the computation of the uncertainty product of zonal $\mathcal S^n$--functions, a generalization of Fern\'andez' ideas from~\cite{nF}.

The paper is organized as follows. After a short background presentation in Section~\ref{sec:sphere}, we introduce polynomial weighted scaling functions and wavelets, present algorithms for reconstruction and decomposition of spherical signals and show the MRA--property in Section~\ref{sec:MRA}. Section~\ref{sec:uncertainty} is devoted to the computation of the uncertainty product of the weighted scaling functions and analysis of its asymptotic behavior. Finally, in Section~\ref{sec:SBA} we discuss the property of positive definiteness of the weighted scaling functions. Two technical lemmas are postponed to the Appendix.

\section{Preliminaries}\label{sec:sphere}

By $\mathcal{S}^n$ we denote the $n$--dimensional unit sphere in $(n+1)$--dimensional Euclidean space~$\mathbb{R}^{n+1}$ with the rotation--invariant measure~$d\sigma$ normalized such that
$$
\Sigma_n=\int_{\mathcal{S}^n}d\sigma=\frac{2\pi^{\lambda+1}}{\Gamma(\lambda+1)},
$$
where~$\lambda$ and~$n$ are related by
$$
\lambda=\frac{n-1}{2}.
$$
The surface element $d\sigma$ is explicitly given by
\begin{equation}\label{eq:surface_element}
d\sigma=\sin^{n-1}\vartheta_1\,\sin^{n-2}\vartheta_2\dots\sin\vartheta_{n-1}d\vartheta_1\,d\vartheta_2\dots d\vartheta_{n-1}d\varphi,
\end{equation}
where $(\vartheta_1,\vartheta_2,\dots,\vartheta_{n-1},\varphi)\in[0,\pi]^{n-1}\times[0,2\pi)$ are spherical coordinates satisfying
\begin{align*}
x_1&=\cos\vartheta_1,\\
x_2&=\sin\vartheta_1\cos\vartheta_2,\\
x_3&=\sin\vartheta_1\sin\vartheta_2\cos\vartheta_3,\\
&\dots\\
x_{n-1}&=\sin\vartheta_1\sin\vartheta_2\dots\sin\vartheta_{n-2}\cos\vartheta_{n-1},\\
x_n&=\sin\vartheta_1\sin\vartheta_2\dots\sin\vartheta_{n-2}\sin\vartheta_{n-1}\cos\varphi,\\
x_{n+1}&=\sin\vartheta_1\sin\vartheta_2\dots\sin\vartheta_{n-2}\sin\vartheta_{n-1}\sin\varphi.
\end{align*}
$\left<x,y\right>$ or $x\cdot y$ stand for the scalar product of vectors with origin in~$O$ and an endpoint on the sphere. As long as it does not lead to misunderstandings, we identify these vectors with points on the sphere.

Scalar product of $f,g\in\mathcal L^2(\mathcal S^n)$ is defined by
$$
\left<f,g\right>_{\mathcal L^2(\mathcal S^n)}=\frac{1}{\Sigma_n}\int_{\mathcal S^n}\overline{f(x)}\,g(x)\,d\sigma(x),
$$
and by $\|\circ\|$ we denote the induced $\mathcal L^2$--norm.

Gegenbauer polynomials $C_l^\lambda$ of order~$\lambda\in\mathbb R$, and degree~$l\in\mathbb{N}_0$, are defined in terms of their generating function
$$
\sum_{l=0}^\infty C_l^\lambda(t)\,r^l=\frac{1}{(1-2tr+r^2)^\lambda},\qquad t\in[-1,1].
$$
They are real-valued and for some fixed $\lambda\ne0$ orthogonal to each other with respect to the weight function~$\left(1-\circ^2\right)^{\lambda-\frac{1}{2}}$, compare~\cite[formula~8.939.8]{GR}.

Let $Q_l$ denote a polynomial on~$\mathbb{R}^{n+1}$ homogeneous of degree~$l$, i.e., such that $Q_l(az)=a^lQ_l(z)$ for all $a\in\mathbb R$ and $z\in\mathbb R^{n+1}$, and harmonic in~$\mathbb{R}^{n+1}$, i.e., satisfying $\Delta Q_l(z)=0$, then $Y_l(x)=Q_l(x)$, $x\in\mathcal S^n$, is called a hyperspherical harmonic of degree~$l$. The set of hyperspherical harmonics of degree~$l$ restricted to~$\mathcal S^n$ is denoted by $\mathcal H_l=\mathcal H_l(\mathcal S^n)$, and we set $\Pi_m=\bigcup_{l=0}^m\mathcal H_l$. $\mathcal H_l$--functions are eigenfunctions of Laplace--Beltrami operator $\Delta^\ast:=\left.\Delta\right|_{\mathcal S^n}$ with eigenvalue $-l(l+2\lambda)=-l(l+n-1)$, further, hyperspherical harmonics of distinct degrees are orthogonal to each other. The number of linearly independent hyperspherical harmonics of degree~$l$ is equal to
$$
N=N(n,l)=\frac{(2l+n-1)(l+n-2)!}{(n-1)!\,l!},
$$
and consequently
$$
\dim\Pi_m=\frac{(n+2m)(n+m-1)!}{n!\,m!}=\frac{n+2m}{n+m}\,\binom{n+m}{m}.
$$
We were not able to find the last result in the literature, therefore, a proof is given in Appendix (Lemma~\ref{lem:dim_Pi}).

Addition theorem states that
\begin{equation}\label{eq:addition_theorem}
C_l^\lambda(x\cdot y)=\frac{\lambda}{l+\lambda}\,\sum_{\kappa=1}^N\overline{Y_l^\kappa(x)}\,Y_l^\kappa(y),
\end{equation}
for any orthonormal set $\{Y_l^\kappa\}_{\kappa=1,2,\dots,N(n,l)}$ of hyperspherical harmonics of degree~$l$ on~$\mathcal S^n$.
In this paper, we will be working with the orthogonal basis for~$\mathcal L^2(\mathcal S^n)=\overline{\bigcup_{l=0}^\infty\mathcal H_l}$, consisting of hyperspherical harmonics given by
\begin{equation}\label{eq:spherical_harmonics}
Y_l^k(x)=A_l^k\prod_{\nu=1}^{n-1}C_{k_{\nu-1}-k_\nu}^{\frac{n-\nu}{2}+k_\nu}(\cos\vartheta_\nu)\sin^{k_\nu}\vartheta_\nu\cdot e^{\pm ik_{n-1}\varphi}
\end{equation}
with $l=k_0\geq k_1\geq\dots\geq k_{n-1}\geq0$, $k$ being a sequence $(k_1,\dots,\pm k_{n-1})$ of integer numbers, and normalization constants~$A_l^k$, compare~\cite{Vilenkin,AI}.

Every $\mathcal{L}^1(\mathcal S^n)$--function~$f$ can be expanded into Laplace series of hyperspherical harmonics by
$$
f\sim\sum_{l=0}^\infty Y_l(f;x),
$$
where $Y_l(f;x)$ is given by
\begin{equation}\label{eq:reprod_kernel_sp}
Y_l(f;x)=\frac{\Gamma(\lambda)(l+\lambda)}{2\pi^{\lambda+1}}\int_{\mathcal S^n}C_l^\lambda(x\cdot y)\,f(y)\,d\sigma(y)
   =\frac{l+\lambda}{\lambda}\left<C_l^\lambda(x\cdot\circ),f\right>.
\end{equation}
For zonal functions (i.e., those depending only on~$\vartheta_1\left<\hat e,x\right>$, where~$\hat e$ is the north-pole of the sphere
$\hat e=(1,0,\dots,0)$) we obtain the Gegenbauer expansion
$$
f(\cos\vartheta_1)=\sum_{l=0}^\infty\widehat f(l)\,C_l^\lambda(\cos\vartheta_1)
$$
with Gegenbauer coefficients
\begin{equation}\label{eq:Gegenbauer_coeffs}
\widehat f(l)=c(l,\lambda)\int_{-1}^1 f(t)\,C_l^\lambda(t)\left(1-t^2\right)^{\lambda-1/2}dt,
\end{equation}
where~$c$ is a constant that depends on~$l$ and~$\lambda$.

For $f,g\in\mathcal L^1(\mathcal S^n)$, $g$ zonal, their convolution $f\ast g$ is defined by
\begin{equation*}
(f\ast g)(x)=\frac{1}{\Sigma_n}\int_{\mathcal S^n}f(y)\,g(x\cdot y)\,d\sigma(y).
\end{equation*}
With this notation we have
\begin{equation}\label{eq:Ylconv}
Y_l(f;x)=\frac{l+\lambda}{\lambda}\,\left(f\ast C_l^\lambda\right)(x),
\end{equation}
i.e., the function $\frac{l+\lambda}{\lambda}\,C_l^\lambda$ is the reproducing kernel for~$\mathcal H_l$.

Further, any function $f\in\mathcal L^2(\mathcal S^n)$ has a unique representation as a mean--convergent series
$$
f(x)=\sum_l\sum_k a_l^kY_l^k(x),\qquad x\in\mathcal S^n,
$$
where
\begin{equation}\label{eq:fourier_coefficient}
a_l^k=a_l^k(f)=\frac{1}{\Sigma_n}\int_{\mathcal S^n}\overline{Y_l^k(x)}\,f(x)\,d\sigma(x)=\left<Y_l^k,f\right>,
\end{equation}
for proof cf.~\cite{Vilenkin}. In analogy to the two-dimensional case, we call $a_l^k$ the Fourier coefficients of the function~$f$.

We identify zonal functions with functions over the interval $[-1,1]$, i.e., whenever it does not lead to mistakes, we write
$$
f(x)=f(\cos\vartheta_1).
$$

\section{Polynomial wavelets and polynomial multiscale decomposition}\label{sec:MRA}

\subsection{The polynomial multiscale decomposition}

In~\cite{MRes} polynomial approximation and polynomial multiscale decomposition of~$\mathcal L^2(\mathcal S^2)$ are discussed. Similarly as in~\cite[Theorem~7]{MRes}, we can argue that for sampling spaces~$\mathcal V_j$, \mbox{$j\in\mathbb{N}$}, given by $\mathcal V_j:=\Pi_{m_j-1}$, where $(m_j)_{j\in\mathbb{N}}$ is a strictly monotonously increasing sequence of positive integers, the family $\{\mathcal V_j\}_{j=1}^\infty$ is a multiscale decomposition of~$\mathcal L^2(\mathcal S^n)$. For our further investigations we set
$$
m_j=2^{j-1},
$$
which yields
$$
\mathcal{V}_j=\Pi_{2^{j-1}-1}=\bigoplus_{l=0}^{2^{j-1}-1}\mathcal H_l
$$
for $j\in\mathbb{N}$. The dimension of~$\mathcal V_j$ is given by
$$
\dim\mathcal V_j=\frac{(n+2^j-2)(n+2^{j-1}-2)!}{n!\,(2^{j-1}-1)!}.
$$
Further, let
$$
I_j:=\left\{(l,k_1,\dots,k_{n-1}):\,l=0,1,\dots,2^{j-1}-1,\,l\geq k_1\geq\dots\geq k_{n-1}\geq0\right\}
$$
be the index set of hyperspherical harmonics which span~$\mathcal V_j$, and
\begin{align*}
\mathcal N_j=\{(s_1,s_2,\dots,s_{n-1},t):\,&s_\nu=0,\dots,2^j\,\text{ for }\,\nu=1,\dots,n-1;\\
  &t=0,\dots,2^{j+1}-1\}.
\end{align*}

We shall prove a sampling theorem analogous to the one given in~\cite{Kernels}. However, since the surface element~$\sigma$ contains higher powers of~$\sin\vartheta_\nu$, compare~\eqref{eq:surface_element}, we need a generalized quadrature formula.

\begin{lem}\label{lem:quadrature} Let~$f$ be a polynomial of degree at most~$M$ and $\alpha\in\mathbb N$. Then
\begin{equation}\label{eq:quadrature}
\int_0^\pi f(\cos\vartheta)\sin^\alpha\vartheta\,d\vartheta=\sum_{u=0}^M\chi_u f\left(\cos\frac{u\pi}{M}\right),
\end{equation}
where
$
\chi_u=\epsilon_u\cdot\omega_u
$
with
$$
\epsilon_0=\epsilon_M=\frac{1}{2}\quad\text{and}\quad\epsilon_u=1\text{ for }u=1,\dots,M-1,
$$
and
$$
\omega_u=\frac{\pi\,\alpha!}{2^{\alpha-1}M}\sum_{\mu=0}^{[M/2]}\frac{\epsilon_{2\mu}(-1)^{\mu}\cos\frac{2\mu u\pi}{M}}{\Gamma\left(\frac{\alpha}{2}-\mu+1\right)\Gamma\left(\frac{\alpha}{2}+\mu+1\right)}.
$$
If $\frac{\alpha}{2}-\mu+1\in-\mathbb N$, the summand is understood to be equal to~$0$.
\end{lem}

\begin{bfseries}Proof.\end{bfseries}
Since~$f$ is a polynomial, $f(\cos\circ)$ has a finite cosine--series representation
$$
f(\cos\vartheta)=\sum_{\mu=0}^M a_\mu\cos(\mu\vartheta),
$$
which coefficients~$a_\mu$ can be computed exactly by the type-I discrete cosine transform with~$M+1$ nodes
$$
a_\mu=\frac{2\epsilon_\mu}{M}\sum_{u=0}^M\epsilon_uf\left(\cos\frac{u\pi}{M}\right)\cos\frac{\mu u\pi}{M}.
$$
The integral of a single frequency~$\cos(\mu\vartheta)$ is given by
$$
\int_0^\pi\cos(\mu\vartheta)\sin^\alpha\vartheta\, d\vartheta=\frac{\pi\,\alpha!\cos\frac{\mu\pi}{2}}{2^\alpha\,\Gamma\left(\frac{\alpha-\mu+2}{2}\right)\Gamma\left(\frac{\alpha+\mu+2}{2}\right)}.
$$
This expression is to be understood as~$0$ for $\mu>\alpha$ with the same parity, and it is equal to~$0$ for odd~$\mu$. Together we obtain
\begin{equation*}\begin{split}
\int_0^\pi f(\cos\vartheta)\sin^\alpha\vartheta\,d\vartheta&=\frac{\pi\,\alpha!}{2^{\alpha-1}M}\sum_{u=0}^M\epsilon_uf\left(\cos\frac{u\pi}{M}\right)\\
&\cdot\sum_{\mu=0}^{[M/2]}\frac{\epsilon_{2\mu}\,(-1)^{\mu}\cos\frac{2\mu u\pi}{M}}{\Gamma\left(\frac{\alpha}{2}-\mu+1\right)\Gamma\left(\frac{\alpha}{2}+\mu+1\right)}.
\end{split}\end{equation*}
\hfill$\Box$

\begin{bfseries}Remark. \end{bfseries}For~$\alpha=1$ this yields the Clenshaw--Curtis quadrature.

\begin{thm}\label{thm:Fourier_coeffs} Let $f\in\mathcal V_j$ be given. Then we have for \mbox{$(l,k)\in I_j$}
\begin{equation}\label{eq:discretized_f_coeffs}
a_l^k=\frac{\pi}{2^j}\sum_{(s,t)\in\mathcal N_j}\left(\prod_{\nu=1}^{n-1}\chi_{s_\nu}^{(j)}\right)f(x_{s,t}^j)\,\overline{Y_l^k(x_{s,t}^j)}
\end{equation}
with
\begin{equation}\label{eq:xst}
x_{s,t}^j=\left(\frac{s_1\pi}{2^j},\dots,\frac{s_{n-1}\pi}{2^j},\frac{t\pi}{2^j}\right)
\end{equation}
and weight functions~$\chi_{s_\nu}^{(j)}$ as given in Lemma~\ref{lem:quadrature} with $M=2^j$ and $\omega_{s_\nu}$ computed for $\alpha=n-\nu$.
\end{thm}

The proof is analogous to the proof of sampling theorem for the two-dimensional case presented in~\cite{Kernels}.

\begin{bfseries}Proof. \end{bfseries}By definition of~$\mathcal V_j$, it suffices to consider the functions \mbox{$f=Y_l^k$, $(l,k)\in I_j$}. Their Fourier coefficients are given by
\begin{equation}\label{eq:fourier_coeffs}\begin{split}
a_{l^\prime}^{k^\prime}=\frac{A_{l^\prime}^{k^\prime}A_l^k}{\Sigma_n}\,
  \prod_{\nu=1}^{n-1}&\int_0^\pi C_{k_{\nu-1}^\prime-k_\nu^\prime}^{\frac{n-\nu}{2}+k_\nu^\prime}(\cos\vartheta_\nu)\,
      C_{k_{\nu-1}-k_\nu}^{\frac{n-\nu}{2}+k_\nu}(\cos\vartheta_\nu)\sin^{k_\nu^\prime+k_\nu+n-\nu}\vartheta_\nu\,d\vartheta_\nu\\
&\cdot\int_0^{2\pi}e^{i(\mp k_{n-1}^\prime\pm k_{n-1})\varphi}\,d\varphi,
\end{split}\end{equation}
compare~\eqref{eq:surface_element}, \eqref{eq:spherical_harmonics}, and~\eqref{eq:fourier_coefficient}.
Since $|\mp k_{n-1}^\prime\pm k_{n-1}|$ is at most equal to \mbox{$2^j-2$}, we have
\begin{align}
\delta_{\pm k_{n-1}^\prime,\pm k_{n-1}}&=\frac{1}{2\pi}\int_0^{2\pi}e^{i(\mp k_{n-1}^\prime\pm k_{n-1})\varphi}\,d\varphi\notag\\
&=\frac{1}{2^{j+1}}\sum_{t=0}^{2^{j+1}-1}\,e^{i(\mp k_{n-1}^\prime\pm k_{n-1})t\pi/2^j}.\label{eq:discrete_e}
\end{align}
Hence, for $k_{n-1}^\prime\ne k_{n-1}$ the theorem is proven.

Now, let~$\nu\in\{1,2,\dots,n-1\}$ be fixed and suppose $k_\nu^\prime=k_\nu$. The product
\begin{equation}\label{eq:product_C}
C_{k_{\nu-1}^\prime-k_\nu}^{\frac{n-\nu}{2}+k_\nu}(\cos\vartheta_\nu)\,C_{k_{\nu-1}-k_\nu}^{\frac{n-\nu}{2}+k_\nu}(\cos\vartheta_\nu)\,\sin^{2k_\nu}\vartheta_\nu
\end{equation}
is a polynomial in~$\cos\vartheta_\nu$ of degree at most~$k_{\nu-1}^\prime+k_{\nu-1}$, which is less than or equal to~$2^j-2$. Therefore, quadrature formula~\eqref{eq:quadrature} with $M=2^j$ and $\alpha=n-\nu$ can be applied and the integral of~\eqref{eq:product_C} with weight $\sin^{n-\nu}\vartheta_\nu$ is equal to
\begin{equation}\label{eq:discrete_C}
\sum_{s_\nu=0}^{2^j}\chi_{s_\nu}\,
   C_{k_{\nu-1}^\prime-k_\nu}^{\frac{n-\nu}{2}+k_\nu}\left(\cos\frac{s_\nu\pi}{2^j}\right)\,
   C_{k_{\nu-1}-k_\nu}^{\frac{n-\nu}{2}+k_\nu}\left(\cos\frac{s_\nu\pi}{2^j}\right)\,\sin^{2k_\nu}\frac{s_\nu\pi}{2^j}.
\end{equation}
Since for~$\lambda\ne0$ the functions $\mapsto C_l^\lambda(\cos\circ)$ and $\mapsto C_{l^\prime}^\lambda(\cos\circ)$ are orthogonal to each other with respect to the weight~$\sin^{2\lambda}\circ$, the sum~\eqref{eq:discrete_C} vanishes for $k_{\nu-1}^\prime\ne k_{\nu-1}$. Thus, by induction, $a_{l^\prime}^{k^\prime}$ vanishes for $(l^\prime,k^\prime)\ne(l,k)$, and if the indices match, the integrals in~\eqref{eq:fourier_coeffs} can be replaced by the sums~\eqref{eq:discrete_e} and~\eqref{eq:discrete_C}. This yields~\eqref{eq:discretized_f_coeffs}.\hfill$\Box$\\

Now we introduce \emph{weighted scaling functions} from~$\mathcal {V}_j$,
\begin{equation}\label{eq:wsf}
\varphi_j(\circ\cdot x_{s,t}^j):=\frac{1}{\sqrt2^{\,nj}}\sum_{(l,k)\in I_j}\overline{Y_l^k(x_{s,t}^j)}\,Y_l^k.
\end{equation}
The addition theorem~\eqref{eq:addition_theorem} yields
\begin{equation}\label{eq:wsf_kernels}
\varphi_j=\frac{1}{\sqrt2^{\,nj}}\sum_{l=0}^{2^{j-1}-1}\frac{l+\lambda}{\lambda}\,C_l^\lambda.
\end{equation}

The next theorem summarizes some of the properties of~$\varphi_j$.

\begin{thm}\label{thm:wsf_properties}Let $(s,t)\in\mathcal N_j$ and~$x_{s,t}^j$ be given by~\eqref{eq:xst}. Then the following holds:
\begin{enumerate}
\item The functions $\varphi_j(\circ\cdot x_{s,t}^j)$ are real--valued.
\item The functions $\varphi_j(\circ\cdot x_{s,t}^j)$ have the reproducing property
\begin{equation}\label{eq:reproducing_property}
\left<\varphi_j(\circ\cdot x_{s,t}^j),f\right>=\frac{1}{\sqrt2^{\,nj}}\,f(x_{s,t}^j)\qquad\text{for all }f\in\mathcal V_j.
\end{equation}
\item We have $\|\varphi_j(\circ\cdot x_{s,t}^j)\|^2=\frac{1}{\sqrt2^{\,nj}}\,\varphi_j(1)=\frac{(n+2^j-2)(n+2^{j-1}-2)!}{2^{nj}\,n!\,(2^{j-1}-1)!}$.
\item The function $\varphi_j(\circ\cdot x_{s,t}^j)$ is localized around~$x_{s,t}^j$, i.e.,
$$
\frac{\|\varphi_j(\circ\cdot x_{s,t}^j)\|}{\varphi_j(x_{s,t}^j\cdot x_{s,t}^j)}=\min\left\{\|f\|:\,f\in\mathcal V_j,\,f(x_{s,t}^j)=1\right\}.
$$
\item It holds
$$
\text{span}\left\{\varphi_j(\circ\cdot x_{s,t}^j):\,(s,t)\in\mathcal N_j\right\}=\mathcal V_j,
$$
and
\begin{equation}\label{eq:discretized_f}
f=\sqrt2^{\,j(n-2)}\pi\sum_{(s,t)\in\mathcal N_j}\left(\prod_{\iota=1}^{n-1}\chi_{s_\iota}^{(j)}\right)f(x_{s,t}^j)\,\varphi_j(\circ\cdot x_{s,t}^j)
\end{equation}
for $f\in\mathcal V_j$.
\item The set
$$
\left\{\left(\prod_{\iota=1}^{n-1}\chi_{s_\iota}^{(j)}\right)^{\!\!1/2}\varphi_j(\circ\cdot x_{s,t}^j):\,(s,t)\in\mathcal N_j\right\}
$$
is a tight frame in~$\mathcal V_j$,
$$
2^{nj-1}\pi\sum_{(s,t)\in\mathcal N_j}\left(\prod_{\iota=1}^{n-1}\chi_{s_\iota}^{(j)}\right)\left|\left<f,\varphi_j(\circ\cdot x_{s,t}^j)\right>\right|^2=\|f\|^2
$$
for every $f\in\mathcal V_j$.
\item The relation
$$
a_l^k\left(\varphi_j(\circ\cdot x_{s,t}^j)\right)=\begin{cases}\frac{1}{\sqrt2^{\,nj}}\,\overline{Y_l^k(x_{s,t}^j)}&\text{for }(l,k)\in I_j,\\0&\text{otherwise}\end{cases}
$$
and the two-scale relation
$$
a_l^k\left(\varphi_j(\circ\cdot x_{s,t}^j)\right)=\begin{cases}\sqrt2^{\,n}a_l^k\left(\varphi_{j+1}(\circ\cdot x_{s,t}^j)\right)&\text{for }(l,k)\in I_j,\\0&\text{otherwise}\end{cases}
$$
hold for the Fourier coefficients of~$\varphi_j$.
\item We have
$$
\int_{\mathcal S^n}\varphi_j(x\cdot x_{s,t}^j)\,d\sigma(x)
   =\frac{2^{3-n(1+\frac{j}{2})}\pi^{\frac{n}{2}+1}\Gamma(2\lambda)}{\Gamma(\lambda)\,\Gamma\!\left(\lambda+\frac{1}{2}\right)\Gamma(\lambda+1)}.
$$
\end{enumerate}
\end{thm}

\begin{bfseries}Remark.\end{bfseries} It follows from property~$3.$ that for big~$j$ the functions $\varphi_j$'s have approximately equal $\mathcal L^2$--norm, more exactly,
$$
\lim_{j\to\infty}\|\varphi_j\|^2=\frac{2^{1-n}}{n!}.
$$

\begin{bfseries}Proof. \end{bfseries}
\begin{enumerate}
\item Follows directly from the representation~\eqref{eq:wsf_kernels}.
\item Use the expression~\eqref{eq:wsf_kernels} for~$\varphi_j$ and the property~\eqref{eq:reprod_kernel_sp} of the functions~$\frac{l+\lambda}{\lambda}\,C_l^\lambda$.
\item It follows from property~2. that
$$
\|\varphi_j(\circ\cdot x_{s,t}^j)\|^2=\frac{1}{\sqrt2^{\,nj}}\,\varphi_j( x_{s,t}^j\cdot x_{s,t}^j)
   =\frac{1}{2^{nj}}\,\sum_{l=0}^{2^{j-1}-1}\frac{l+\lambda}{\lambda}\,C_l^\lambda(1),
$$
and further, by~\cite[formula~8.937.4]{GR},
\begin{align*}
\|\varphi_j(\circ\cdot x_{s,t}^j)\|^2&=\frac{1}{2^{nj}}\sum_{l=0}^{2^{j-1}-1}\frac{l+\lambda}{\lambda}\,\binom{l+2\lambda-1}{l}\\
&=\frac{1}{2^{nj}}\sum_{l=0}^{2^{j-1}-1}\frac{(2l+n-1)(l+n-2)!}{(n-1)!\,l!}.
\end{align*}
Consequently, by Lemma~\ref{lem:dim_Pi}, we have
$$
\|\varphi_j(\circ\cdot x_{s,t}^j)\|^2=\frac{(n+2^j-2)(n+2^{j-1}-2)!}{2^{nj}\,n!\,(2^{j-1}-1)!}.
$$
\item For all $f\in\mathcal V_j$ with $f(x_{s,t}^j)=1$ we have
$$
1=\sum_{(l,k)\in I_j}a_l^k(f)\,Y_l^k(x_{s,t}^j).
$$
The Cauchy--Schwarz inequality yields
\begin{equation}\label{eq:CS}
1\leq\sum_{(l,k)\in I_j}|a_l^k(f)|^2\cdot\sum_{(l,k)\in I_j}|Y_l^k(x_{s,t}^j)|^2,
\end{equation}
and equality in~\eqref{eq:CS} holds only for a function~$\tilde f$ such that $\left(a_l^k(\tilde f)\right)_{(l,k)\in I_j}$ and $\left(\overline{Y_l^k(x_{s,t}^j)}\right)_{(l,k)\in I_j}$ are linearly dependent, i.e.,
$$
a_l^k(\tilde f)=\alpha\,\overline{Y_l^k(x_{s,t}^j})
$$
for all $(l,k)\in I_j$ and for some constant $\alpha\in\mathbb C$. Consequently, $\tilde f$ is an $\frac{\sqrt2^{\,nj}\alpha}{\Sigma_n}$--multiple of~$\varphi_j$, and for all $f\in\mathcal V_j$ with $f(x_{s,t}^j)=1$ we have
$$
\|f\|^2\geq\|\tilde f\|^2=\left\|\frac{\varphi_j(\circ\cdot x_{s,t}^j)}{\varphi_j(1)}\right\|.
$$
\item Let~$f\in\mathcal V_j$ be given by its Fourier series,
$$
f=\sum_{(l,k)\in I_j}a_l^k\,Y_l^k.
$$
The coefficients~$a_l^k$ can be computed by~\eqref{eq:discretized_f_coeffs}. Together with~\eqref{eq:wsf} this yields~\eqref{eq:discretized_f}, i.e., $\mathcal V_j$ is spanned by the set $\left\{\varphi_j(\circ\cdot x_{s,t}^j):\,(s,t)\in\mathcal N_j\right\}$.
\item For the $\mathcal L^2$-- norm of~$f\in\mathcal V_j$ we write
$$
\|f\|^2=\left<f,f\right>=\sum_{(l,k)\in I_j}\overline{a_l^k}\left<Y_l^k,f\right>.
$$
From~\eqref{eq:discretized_f} it follows
$$
\|f\|^2=2^{\frac{j(n-2)}{2}}\pi\sum_{(l,k)\in I_j}\sum_{(s,t)\in\mathcal N_j}\overline{a_l^k}\left(\prod_{\iota=1}^{n-1}\chi_{s_\iota}^{(j)}\right)
   f(x_{s,t}^j)\left<Y_l^k,\varphi_j(\circ\cdot x_{s,t}^j)\right>.
$$
Further, we apply~\eqref{eq:reproducing_property} for~$f(x_{s,t})$ and for $\left<Y_l^k,\varphi_j(\circ\cdot x_{s,t}^j)\right>$, and we obtain
\begin{align*}
\|f\|^2&=2^{\frac{j(n-2)}{2}}\pi\sum_{(l,k)\in I_j}\sum_{(s,t)\in\mathcal N_j}\overline{a_l^k}\left(\prod_{\iota=1}^{n-1}\chi_{s_\iota}^{(j)}\right)
   \left<\varphi_j(\circ\cdot x_{s,t}^j),f\right>\overline{Y_l^k(x_{s,t}^j)}\\
&=2^{\frac{j(n-2)}{2}}\pi\sum_{(s,t)\in\mathcal N_j}\left(\prod_{\iota=1}^{n-1}\chi_{s_\iota}^{(j)}\right) \left<\varphi_j(\circ\cdot x_{s,t}^j),f\right>
   \sum_{(l,k)\in I_j}\overline{a_l^kY_l^k(x_{s,t}^j)}.
\end{align*}
The last sum in this expression is equal to~$\overline{f(x_{s,t}^j)}$, and using for it again the reproducing property~\eqref{eq:reproducing_property} we obtain
$$
\|f\|^2=2^{nj-1}\pi\sum_{(s,t)\in\mathcal N_j}\left(\prod_{\iota=1}^{n-1}\chi_{s_\iota}^{(j)}\right)\left|\left<\varphi_j(\circ\cdot x_{s,t}^j),f\right>\right|^2.
$$
\item Follows directly from the representation~\eqref{eq:wsf}.
\item The integral over~$\mathcal S^n$ is invariant with respect to rotations, therefore, we can fix $x_{s,t}^j=\hat e$. Thus,
$$
\int_{\mathcal S^n}\varphi_j(x\cdot x_{s,t}^j)\,d\sigma(x)=\Sigma_{n-1}\int_0^\pi\varphi_j(\cos\vartheta_1)\sin^{n-1}\vartheta_1d\vartheta_1.
$$
From~\eqref{eq:wsf_kernels} we conclude
$$
\int_{\mathcal S^n}\varphi_j(x\cdot x_{s,t}^j)\,d\sigma(x)
   =\frac{\Sigma_{n-1}}{\sqrt 2^{\,nj}}\sum_{l=0}^{2^{j-1}-1}\frac{l+\lambda}{\lambda}\int_0^\pi C_l^\lambda(\cos\vartheta_1)\sin^{n-1}\vartheta_1d\vartheta_1,
$$
\end{enumerate}
and further, since $C_0^\lambda(t)=1$ (compare~\cite[formula~8.937.3]{GR}),
\begin{align*}
\int_{\mathcal S^n}\varphi_j(x\cdot x_{s,t}^j)\,d\sigma(x)
   &=\frac{\Sigma_{n-1}}{\sqrt 2^{\,nj}}
   \sum_{l=0}^{2^{j-1}-1}\frac{l+\lambda}{\lambda}\int_{-1}^1C_0^\lambda(t)\,C_l^\lambda(t)\left(1-t^2\right)^{\lambda-\frac{1}{2}}dt\\
&=\frac{2\pi^{\lambda+\frac{1}{2}}}{2^{\left(\lambda+\frac{1}{2}\right)j}\,\Gamma\left(\lambda+\frac{1}{2}\right)}
   \sum_{l=0}^{2^{j-1}-1}\frac{l+\lambda}{\lambda}\cdot\delta_{l,0}\cdot\frac{2^{1-2\lambda}\pi\,\Gamma(2\lambda)}{\Gamma(\lambda)\Gamma(\lambda+1)}.
\end{align*}
The last equality holds by~\cite[formula~8.939.6]{GR}. Consequently,
$$
\int_{\mathcal S^n}\varphi_j(x\cdot x_{s,t}^j)\,d\sigma(x)
   =\frac{2^{2-2\lambda-j\lambda-\frac{1}{2}j}\pi^{\lambda+\frac{3}{2}}\Gamma(2\lambda)}{\Gamma(\lambda)\,\Gamma\!\left(\lambda+\frac{1}{2}\right)\Gamma(\lambda+1)}.
$$
\hfill$\Box$

\subsection{The polynomial wavelet space}
Wavelet spaces~$\mathcal W_j$, $j\in\mathbb N$, are defined as the direct sum
$$
\mathcal W_j=\bigoplus_{l=2^{j-1}}^{2^j-1}\mathcal H_l.
$$
Their dimensions are given by
$$
\dim\mathcal W_j=\frac{1}{n!}\left(\frac{(n+2^{j+1}-2)(n+2^j-2)!}{(2^j-1)!}-\frac{(n+2^j-2)(n+2^{j-1}-2)!}{(2^{j-1}-1)!}\right).
$$
We introduce the \emph{weighted wavelets} from~$\mathcal W_j$,
\begin{equation}\label{eq:definition_wavelets}
\psi_j(\circ\cdot x_{s,t}^{j+1}):=\frac{1}{\sqrt 2^{\,nj}}\sum_{(l,k)\in I_{j+1}\setminus I_j}\overline{Y_l^k(x_{s,t}^{j+1})}\,Y_l^k,\qquad(s,t)\in\mathcal N_{j+1}.
\end{equation}
Consequently, by~\eqref{eq:addition_theorem} we have
\begin{equation}\label{eq:wavelets_kernels}
\psi_j=\frac{1}{\sqrt 2^{\,nj}}\sum_{l=2^{j-1}}^{2^j-1}\frac{l+\lambda}{\lambda}\,C_l^\lambda.\\
\end{equation}

Their properties are summarized in the next theorem.

\begin{thm}\label{thm:ww_proprties}Let $(s,t)\in\mathcal N_{j+1}$ and~$x_{s,t}^{j+1}$ be given by~\eqref{eq:xst}. Then the following holds:
\begin{enumerate}
\item The functions $\psi_j(\circ\cdot x_{s,t}^{j+1})$ are real--valued.
\item The functions $\psi_j(\circ\cdot x_{s,t}^{j+1})$ have the reproducing property
\begin{equation}\label{eq:wv_reproducing_property}
\left<\psi_j(\circ\cdot x_{s,t}^{j+1}),f\right>=\frac{1}{\sqrt 2^{\,nj}}f(x_{s,t}^{j+1})\qquad\text{for all }f\in\mathcal W_j.
\end{equation}
\item $\varphi_j$ and $\psi_j$ are orthogonal to each other,
$$
\left<\varphi_j(x_{s^\prime,t^\prime}^j\cdot\circ),\psi_j(x_{s,t}^{j+1}\cdot\circ)\right>=0\qquad\text{for all }(s^\prime,t^\prime)\in\mathcal N_j,\,(s,t)\in\mathcal N_{j+1}.
$$
\item We have
\begin{align*}
&\|\psi_j(\circ\cdot x_{s,t}^{j+1})\|^2=\frac{1}{\sqrt 2^{\,nj}}\,\psi_j(1)\\
&\qquad=\frac{1}{2^{nj}n!}\left(\frac{(n+2^{j+1}-2)(n+2^j-2)!}{2^n\,(2^j-1)!}-\frac{(n+2^j-2)(n+2^{j-1}-2)!}{(2^{j-1}-1)!}\right).
\end{align*}
\item The function $\psi_j(\circ\cdot x_{s,t}^{j+1})$ is localized around~$x_{s,t}^{j+1}$, i.e.,
$$
\frac{\|\psi_j(\circ\cdot x_{s,t}^{j+1})\|}{\psi_j(x_{s,t}^{j+1}\cdot x_{s,t}^{j+1})}=\min\left\{\|f\|:\,f\in\mathcal W_j,\,f(x_{s,t}^{j+1})=1\right\}.
$$
\item It holds
$$
\text{span}\left\{\psi_j(\circ\cdot x_{s,t}^{j+1}):\,(s,t)\in\mathcal N_{j+1}\right\}=\mathcal W_j,
$$
and
\begin{equation}\label{eq:discretized_finW}
f=\sqrt2^{\,(n-2)j-2}\,\pi\sum_{(s,t)\in\mathcal N_{j+1}}\left(\prod_{\iota=1}^{n-1}\chi_{s_\iota}^{(j+1)}\right)f(x_{s,t}^{j+1})\,\psi_j(\circ\cdot x_{s,t}^{j+1})
\end{equation}
for $f\in\mathcal W_j$.
\item The set
$$
\left\{\left(\prod_{\iota=1}^{n-1}\chi_{s_\iota}^{(j)}\right)^{\!\!1/2}\psi_j(\circ\cdot x_{s,t}^{j+1}):\,(s,t)\in\mathcal N_{j+1}\right\}
$$
is a tight frame in~$\mathcal W_j$,
$$
\sqrt2^{\,(n-2)j-2}\,\pi\sum_{(s,t)\in\mathcal N_{j+1}}\left(\prod_{\iota=1}^{n-1}\chi_{s_\iota}^{(j+1)}\right)\left|\left<f,\psi_j(\circ\cdot x_{s,t}^{j+1})\right>\right|^2=\|f\|^2
$$
for every $f\in\mathcal W_j$.
\item The relation
$$
a_l^k\left(\psi_j(\circ\cdot x_{s,t}^{j+1})\right)=\begin{cases}\frac{1}{\sqrt 2^{\,nj}}\,\overline{Y_l^k(x_{s,t}^{j+1})}&\text{for }(l,k)\in I_{j+1}\setminus I_j,\\0&\text{otherwise}\end{cases}
$$
and the two-scale realtion
$$
a_l^k\left(\psi_j(\circ\cdot x_{s,t}^{j+1})\right)=\begin{cases}a_l^k\left(\varphi_{j+1}(\circ\cdot x_{s,t}^{j+1})\right)&\text{for }(l,k)\in I_{j+1}\setminus I_j,\\0&\text{otherwise}\end{cases}
$$
hold for the Fourier coefficients of~$\psi_j$.
\item We have
$$
\int_{\mathcal S^n}\psi_j(x\cdot x_{s,t}^{j+1})\,d\sigma(x)=0.
$$
\end{enumerate}
\end{thm}

\begin{bfseries}Proof. \end{bfseries}Item 3. follows directly from the definitions of~$\varphi_j$ and~$\psi_j$. Proof of the other items is analogous to the proof of Theorem~\ref{thm:wsf_properties}.

\begin{bfseries}Remark. \end{bfseries}For big values of~$j$, $\psi_j$ vanishes in $\mathcal L^2$--norm, more exactly,
$$
\|\psi_j\|^2=\mathcal O(2^{-j}),\qquad j\to\infty.
$$

\subsection{Algorithms for reconstruction and decomposition}\label{subs:algorithms}
Let $v_{j+1}\in\mathcal V_{j+1}$, $v_j\in\mathcal V_j$, and $w_j\in\mathcal W_j$ be given. According to~\eqref{eq:discretized_f} and~\eqref{eq:discretized_finW}, these functions can be represented by
\begin{align*}
v_{j+1}&=\sqrt2^{\,(j+1)(n-2)}\,\pi
   \sum_{(s,t)\in\mathcal N_{j+1}}\left(\prod_{\iota=1}^{n-1}\chi_{s_\iota}^{(j+1)}\right)v_{j+1}(x_{s,t}^{j+1})\,\varphi_{j+1}(\circ\cdot x_{s,t}^{j+1}),\\
v_j&=\sqrt2^{\,j(n-2)}\,\pi\sum_{(s,t)\in\mathcal N_j}\left(\prod_{\iota=1}^{n-1}\chi_{s_\iota}^{(j)}\right)v_j(x_{s,t}^j)\,\varphi_j(\circ\cdot x_{s,t}^j),\\
w_j&=\sqrt2^{\,j(n-2)-2}\pi\sum_{(s,t)\in\mathcal N_{j+1}}\left(\prod_{\iota=1}^{n-1}\chi_{s_\iota}^{(j+1)}\right)w_j(x_{s,t}^{j+1})\,\psi_j(\circ\cdot x_{s,t}^{j+1}).
\end{align*}
Decomposition of the function~$v_{j+1}$ into functions~$v_j$ and~$w_j$,
$$
v_{j+1}=v_j+w_j
$$
is unique and the Fourier coefficients~$a_l^k$ satisfy
$$
a_l^k(v_{j+1})=\begin{cases}a_l^k(v_j)+a_l^k(w_j)&\text{for }(l,k)\in I_{j+1},\\0&\text{for }(l,k)\notin I_{j+1}.\end{cases}
$$
Consequently,
\begin{align*}
a_l^k(v_j)=a_l^k(v_{j+1})\qquad&\text{for }(l,k)\in I_j,\\
a_l^k(w_j)=a_l^k(v_{j+1})\qquad&\text{for }(l,k)\in I_{j+1}\setminus I_j.
\end{align*}

It follows from Theorem~\ref{thm:Fourier_coeffs} that
\begin{align}
v_j(x)&=\sum_{(l,k)\in I_j}a_l^k(v_{j+1})Y_l^k(x)\notag\\
&=\frac{\pi}{2^j}\sum_{\genfrac{}{}{0pt}{}{(l,k)\in I_j}{(s,t)\in\mathcal N_j}}\left(\prod_{\iota=1}^{n-1}\chi_{s_\iota}^{(j)}\right)
   v_{j+1}(x_{s,t}^j)\,\overline{Y_l^k(x_{s,t}^j)}\,Y_l^k(x)\label{eq:mapping_vj}
\end{align}
and
\begin{align}
w_j(x)&=\sum_{(l,k)\in I_{j+1}\setminus I_j}a_l^k(v_{j+1})Y_l^k(x)\notag\\
&=\frac{\pi}{2^{j+1}}\sum_{\genfrac{}{}{0pt}{}{(l,k)\in I_{j+1}\setminus I_j}{(s,t)\in\mathcal N_{j+1}}}
   \left(\prod_{\iota=1}^{n-1}\chi_{s_\iota}^{(j+1)}\right)v_{j+1}(x_{s,t}^{j+1})\,\overline{Y_l^k(x_{s,t}^{j+1})}\,Y_l^k(x)\label{eq:mapping_wj}
\end{align}
for all $x\in\mathcal S^n$, and also particularly for~$x_{s,t}^{(j+1)}$ with $(s,t)\in\mathcal N_{j+1}$. Now let~$R$ be the operator
$$
R\left(\mathbf v^{(j+1)}\right)=\mathbf v^{(j)}
$$
mapping
$$
\mathbf v^{(j+1)}=\left(v_{j+1}(x_{s,t}^{j+1})\right)_{(s,t)\in\mathcal N_{j+1}}
$$
onto
$$
\mathbf v^{(j)}=\left(v_j(x_{s,t}^j)\right)_{(s,t)\in\mathcal N_j}
$$
via~\eqref{eq:mapping_vj}, and analogously, let~$Q$ be the operator
$$
Q\left(\mathbf v^{(j+1)}\right)=\mathbf w^{(j)}
$$
defined by~\eqref{eq:mapping_wj} for
$$
\mathbf w^{(j)}=\left(w_j(x_{s,t}^{j+1})\right)_{(s,t)\in\mathcal N_{j+1}}.
$$
The operators~$R$ and~$Q$ describe the decomposition of a function $v_{j+1}\in\mathcal V_{j+1}$.\\

\newlength{\ii}
\setlength{\ii}{\textwidth}
\noindent
\fbox{
\begin{minipage}{\ii}
\hspace*{10em}\begin{bfseries}Decomposition Algorithm\end{bfseries}\\
Input: $$\mathbf v^{(j+1)}$$
Compute for $i=j,j-1,\dots,1$: $$\mathbf v^{(i)}=R\left(\mathbf v^{(i+1)}\right)$$$$\mathbf w^{(i)}=Q\left(\mathbf v^{(i+1)}\right)$$
Output: $$\mathbf v^{(1)},\qquad\mathbf w^{(i)},\,i=1,2,\dots,j$$
\end{minipage}}\\

The functions~$v_1$ and $w_i$, $i=1,2,\dots,j$, can be reconstructed from their samples $\mathbf v^{(1)}$, $\mathbf w^{(i)}$, $i=1,2,\dots,j$, via~\eqref{eq:mapping_vj} and~\eqref{eq:mapping_wj}.\\

The reconstruction of a function $v_{j+1}\in\mathcal V_{j+1}$
$$
v_j+w_j=v_{j+1}
$$
is given by
\begin{align*}
\mathbf v^{(j+1)}&=R^\ast\left(\mathbf v^{(j)}\right)+Q^\ast\left(\mathbf w^{(j)}\right)
\end{align*}
with
$$
\left(R^\ast\left(\mathbf v^{(j)}\right)\right)_{p,q}=\frac{\pi}{2^j}\sum_{\genfrac{}{}{0pt}{}{(l,k)\in I_j}{(s,t)\in\mathcal N_j}}
   \left(\prod_{\iota=1}^{n-1}\chi_{s_\iota}^{(j)}\right)(\mathbf v^{(j)})_{s,t}\,\overline{Y_l^k(x_{s,t}^j)}\,Y_l^k(x_{p,q}^{j+1})
$$
and
$$
\left(Q^\ast\left(\mathbf w^{(j)}\right)\right)_{p,q}=\frac{\pi}{2^{j+1}}\sum_{\genfrac{}{}{0pt}{}{(l,k)\in I_{j+1}\setminus I_j}{(s,t)\in\mathcal N_{j+1}}}
   \left(\prod_{\iota=1}^{n-1}\chi_{s_\iota}^{(j+1)}\right)(\mathbf w^{(j)})_{s,t}\,\overline{Y_l^k(x_{s,t}^{j+1})}\,Y_l^k(x_{p,q}^{j+1}),
$$
for $(p,q)\in\mathcal N_{j+1}$. The operators $R^\ast$ and $Q^\ast$ are adjoints of~$R$ and~$Q$, respectively.\\

\noindent
\fbox{
\begin{minipage}{\ii}
\hspace*{10em}\begin{bfseries}Reconstruction Algorithm\end{bfseries}\\
Input: $$\mathbf v^{(1)},\qquad\mathbf w^{(i)},\,i=1,2,\dots,j$$
Compute for $i=1,2,\dots,j$: $$\mathbf v^{(i+1)}=R^\ast\left(\mathbf v^{(i)}\right)+Q^\ast\left(\mathbf w^{(i)}\right)$$
Output: $$\mathbf v^{(j+1)}$$
\end{minipage}}\\

\subsubsection{Multiresolution analysis}

The sequence $(\mathcal V_j)_{j=1}^\infty$ is a multiscale decomposition of~$\mathcal L^2(\mathcal S^n)$, i.e.,
\begin{enumerate}
\item $\mathcal V_j\subset\mathcal V_{j+1}$ for $j\in\mathbb N$ and
\item closure$\left(\bigcup_{j\in\mathbb N}\mathcal V_j\right)=\mathcal L^2(\mathcal S^n)$.
\end{enumerate}

The first property results from definition, the second one holds true by Stone--Weierstrass theorem since~$\mathcal S^n$ is compact.

\section{Uncertainty of weighted scaling functions}\label{sec:uncertainty}

In order to give a deepen characterization of $\varphi_j$'s and $\psi_j$'s, we compute their uncertainty product. According to~\cite{RV97} the variances in space and momentum domain of a $\mathcal C^2(\mathcal S^n)$--function~$f$ with $\int_{\mathcal S^n}x\,|f(x)|^2\,d\sigma(x)\ne0$ are given by
$$
\text{var}_S(f)=\left(\frac{\int_{\mathcal S^n}|f(x)|^2\,d\sigma(x)}{\int_{\mathcal S^n}x\,|f(x)|^2\,d\sigma(x)}\right)^2-1
$$
and
$$
\text{var}_M(f)=-\frac{\int_{\mathcal S^n}\Delta^\ast f(x)\cdot \bar f(x)\,d\sigma(x)}{\int_{\mathcal S^n}|f(x)|^2\,d\sigma(x)},
$$
where $\Delta^\ast$ is the Laplace--Beltrami operator on~$\mathcal S^n$, compare also~\cite{GG04} and~\cite{nF}.

\begin{df} The quantity
$$
U(f)=\sqrt{\text{var}_S(f)}\cdot\sqrt{\text{var}_M(f)}
$$
is called uncertainty product of~$f$.
\end{df}

According to \cite[Theorem~1.2]{RV97}, for $SO(n)$--invariant functions $f\in\mathcal L^2(\mathcal S^n)\cap\mathcal C^2(\mathcal S^n)$
$$
U(f)\geq\frac{n}{2}
$$
and the lower bound is optimal (compare also~\cite{wE11} for the non-zonal case). It is the limiting value for $t\to0$ of the uncertatinty product of the so-called Gaussian measures
$$
G_t^\lambda=\frac{\widetilde G_t^\lambda}{\|\widetilde G_t^\lambda\|_1},\qquad
   \widetilde G_t^\lambda=\sum_{l=0}^\infty e^{-\frac{t\,l(l+2\lambda)}{2}}\,\frac{C_l^\lambda}{\|C_l^\lambda\|_2^2},
$$
cf. \cite[Proposition 3.3]{RV97}.

\begin{lem} Let a zonal $\mathcal L^2(\mathcal S^n)$--function be given by its Gegenbauer expansion
$$
f(t)=\sum_{l=0}^\infty\widehat f(l)\,C_l^\lambda(t).
$$
Its variances in space and momentum domain are given by
\begin{align}
\text{var}_S(f)&=\left(\frac{\sum_{l=0}^\infty\frac{\lambda}{l+\lambda}\,\binom{l+2\lambda-1}{l}\,|\widehat f(l)|^2}{\sum_{l=0}^\infty\binom{l+2\lambda}{l}\,
   \frac{\lambda^2\left[\overline{\widehat f(l)}\,\widehat f(l+1)+\widehat f(l)\,\overline{\widehat f(l+1)}\right]}{(l+\lambda)(l+\lambda+1)}}\right)^2-1,\label{eq:varS}\\
\text{var}_M(f)&=\frac{\sum_{l=1}^\infty\frac{ l\lambda(l+2\lambda)}{l+\lambda}\,\binom{l+2\lambda-1}{l}\,|\widehat f(l)|^2}
   {\sum_{l=0}^\infty\frac{\lambda}{l+\lambda}\,\binom{l+2\lambda-1}{l}\,|\widehat f(l)|^2}\label{eq:varM}
\end{align}
whenever the series are convergent.
\end{lem}

\begin{bfseries}Proof. \end{bfseries}Since $\frac{l+\lambda}{\lambda}\,C_l^\lambda$ is the reproducing kernel of~$\mathcal H_l$, we obtain
\begin{align*}
\frac{1}{\Sigma_n}\int_{\mathcal S^n}&|f(x)|^2\,d\sigma(x)=\overline f\ast f(\hat e)
   =\sum_{l=0}^\infty\overline{\widehat f(l)}\,C_l^\lambda\ast\sum_{l=0}^\infty\widehat f(l)\,C_l^\lambda(1)\\
&=\sum_{l=0}^\infty\frac{\lambda}{l+\lambda}\,|\widehat f(l)|^2\,C_l^\lambda(1)
\end{align*}
For $C_l^\lambda(1)$ we write $\binom{l+2\lambda-1}{l}$, compare~\cite[formula~8.937.4]{GR}, and obtain
$$
\frac{1}{\Sigma_n}\int_{\mathcal S^n}|f(x)|^2\,d\sigma(x)=\sum_{l=0}^\infty\frac{\lambda}{l+\lambda}\,\binom{l+2\lambda-1}{l}\,|\widehat f(l)|^2.
$$

In order to compute $\int_{\mathcal S^n}x\,|f(x)|^2\,d\sigma(x)$, denote by~$g$ the function $x\mapsto xf(x)$ and by~$\widetilde C_l^\lambda$ the function~$t\mapsto t\,C_l^\lambda(t)$. We have
$$
\frac{1}{\Sigma_n}\int_{\mathcal S^n}x\,|f(x)|^2\,d\sigma(x)=\overline f\ast g(\hat e)
   =\sum_{l=0}^\infty\overline{\widehat f(l)}\,C_l^\lambda\ast\sum_{l=0}^\infty\widehat f(l)\,\widetilde C_l^\lambda(1).
$$
Now, for calculation of~$\widetilde C_l^\lambda$ we use formulae~8.932.1 and~8.930 from~\cite{GR}:
\begin{align*}
(l+1)\,C_{l+1}^{\lambda}(t)&=2(l+\lambda)t\,C_l^{\lambda}(t)-(l+2\lambda-1)C_{l-1}^{\lambda}(t)\qquad\text{for }l\geq1,\\
C_1^\lambda(t)&=2\lambda t=2\lambda\,t C_0^\lambda(t),
\end{align*}
and obtain
\begin{align*}
\frac{1}{\Sigma_n}\int_{\mathcal S^n}&x\,|f(x)|^2\,d\sigma(x)=\frac{1}{2\lambda}\sum_{l=0}^\infty\overline{\widehat f(l)}\,C_l^\lambda\ast\widehat f(0)\,C_1^\lambda(1)\\
&+\sum_{l=0}^\infty\overline{\widehat f(l)}\,C_l^\lambda\ast\sum_{l=1}^\infty\frac{l+1}{2(l+\lambda)}\,\widehat f(l)\,C_{l+1}^\lambda(1)\\
&+\sum_{l=0}^\infty\overline{\widehat f(l)}\,C_l^\lambda\ast\sum_{l=1}^\infty\frac{l+2\lambda-1}{2(l+\lambda)}\,\widehat f(l)\,C_{l-1}^\lambda(1).
\end{align*}
Further, by the reproducing kernel property of~$\frac{l+\lambda}{\lambda}\,C_l^\lambda$ we get
\begin{align*}
\frac{1}{\Sigma_n}\int_{\mathcal S^n}&x\,|f(x)|^2\,d\sigma(x)=\frac{1}{2(\lambda+1)}\,\overline{\widehat f(1)}\,\widehat f(0)\,C_1^\lambda(1)\\
&+\sum_{l=1}^\infty\frac{\lambda(l+1)}{2(l+\lambda)(l+\lambda+1)}\,\overline{\widehat f(l+1)}\,\widehat f(l)\,C_{l+1}^\lambda(1)\\
&+\sum_{l=1}^\infty\frac{\lambda(l+2\lambda-1)}{2(l+\lambda)(l+\lambda-1)}\,\overline{\widehat f(l-1)}\,\widehat f(l)\,C_{l-1}^\lambda(1),
\end{align*}
and consequently, by index--shift $l\mapsto l-1$ in the second series,
\begin{align*}
\frac{1}{\Sigma_n}\int_{\mathcal S^n}&x\,|f(x)|^2\,d\sigma(x)=\frac{\lambda}{\lambda+1}\,\overline{\widehat f(1)}\,\widehat f(0)\\
&+\sum_{l=1}^\infty\frac{\lambda(l+1)}{2(l+\lambda)(l+\lambda+1)}\,\binom{l+2\lambda}{l+1}\,\overline{\widehat f(l+1)}\,\widehat f(l)\\
&+\sum_{l=0}^\infty\frac{\lambda(l+2\lambda)}{2(l+\lambda+1)(l+\lambda)}\,\binom{l+2\lambda-1}{l}\,\overline{\widehat f(l)}\,\widehat f(l+1).
\end{align*}
Altogether we have
\begin{align*}
\frac{1}{\Sigma_n}\int_{\mathcal S^n}&x\,|f(x)|^2\,d\sigma(x)\\
&=\sum_{l=0}^\infty\frac{\lambda^2}{(l+\lambda)(l+\lambda+1)}\binom{l+2\lambda}{l}\,
   \left[\overline{\widehat f(l+1)}\,\widehat f(l)+\overline{\widehat f(l)}\,\widehat f(l+1)\right].
\end{align*}
The same formula (up to complex conjugation and for zero--mean functions) is derived in~\cite{DX12} with slightly different methods.

Finally, since
$$
-\Delta^\ast f=\sum_{l=1}^\infty l(l+2\lambda)\widehat f(l)\,C_l^\lambda,
$$
we have
\begin{align*}
\frac{-1}{\Sigma_n}\int_{\mathcal S^n}&\Delta^\ast f(x)\cdot \bar f(x)\,d\sigma(x)
   =\sum_{l=1}^\infty l(l+2\lambda)\widehat f(l)\,C_l^\lambda\ast\sum_{l=0}^\infty\overline{\widehat f(l)}\,C_l^\lambda(1)\\
&=\sum_{l=1}^\infty l(l+2\lambda)|\widehat f(l)|^2\cdot\frac{\lambda}{l+\lambda}\,C_l^\lambda(1)\\
&=\sum_{l=1}^\infty\frac{ l\lambda(l+2\lambda)}{l+\lambda}\,\binom{l+2\lambda-1}{l}\,|\widehat f(l)|^2.
\end{align*}
Combining the above formulae we obtain~\eqref{eq:varS} and~\eqref{eq:varM}.\hfill$\Box$

\begin{bfseries}Remark. \end{bfseries}For $\lambda=\frac{1}{2}$ these formulae reduce to those obtained in~\cite[Section~5.2]{nF}.

Now, we want to compute the space and momentum variances for the weighted scaling functions. However, in order to obtain more general results, we do not restrict ourselves to the scales $2^{j-1}-1$. Let us introduce the notation
$$
\Phi_m=C(m)\cdot\sum_{l=0}^m\frac{l+\lambda}{\lambda}\,C_l^\lambda,
$$
where~$C(m)$ is a constant.

\begin{prop} The variances in space and momentum domain of~$\Phi_m$, $m\in\mathbb N$, are equal to
\begin{align}
\text{var}_S(\Phi_m)&=\left(\frac{2m+2\lambda+1}{2m}\right)^2-1,\label{eq:varSPhim}\\
\text{var}_M(\Phi_m)&=\frac{m\,(m+2\lambda+1)\,(2\lambda+1)}{2\lambda+3}.\label{eq:varMPhim}
\end{align}
\end{prop}

\begin{bfseries}Proof. \end{bfseries}In the case of weighted scaling functions we have $\widehat f(l)=\frac{l+\lambda}{\sqrt2^{\,nj}\lambda}$ for $l\leq m=2^{j-1}-1$. With these values of Gegenbauer coefficients we obtain
\begin{align*}
\frac{2^{\,nj}}{\Sigma_n}\int_{\mathcal S^n}&|f(x)|^2\,d\sigma(x)
   =2^{\,nj}\cdot\sum_{l=0}^m\frac{\lambda}{l+\lambda}\,\binom{l+2\lambda-1}{l}\,|\widehat f(l)|^2\\
&=\sum_{l=0}^m\frac{l+\lambda}{\lambda}\,\binom{l+2\lambda-1}{l}
   =\frac{(2m+2\lambda+1)\,(m+2\lambda)!}{m!\,(2\lambda+1)!},
\end{align*}
compare Lemma~\ref{lem:dim_Pi}, and
\begin{align*}
\frac{2^{\,nj}}{\Sigma_n}\int_{\mathcal S^n}&x\,|f(x)|^2\,d\sigma(x)\\
&=2^{\,nj}\cdot\sum_{l=0}^{m-1}\binom{l+2\lambda}{l}\,
   \frac{\lambda^2\left[\overline{\widehat f(l)}\,\widehat f(l+1)+\widehat f(l)\,\overline{\widehat f(l+1)}\right]}{(l+\lambda)(l+\lambda+1)}\\
&=2\cdot\sum_{l=0}^{m-1}\binom{l+2\lambda}{2\lambda}=2\,\binom{m+2\lambda}{2\lambda+1},
\end{align*}
compare~\cite[formula~0.15.1]{GR}, and finally
\begin{align*}
\frac{-2^{\,nj}}{\Sigma_n}\int_{\mathcal S^n}&\Delta^\ast f(x)\cdot \bar f(x)\,d\sigma(x)\\
&=2^{\,nj}\cdot\sum_{l=0}^{m}\,\frac{ l\lambda(l+2\lambda)}{l+\lambda}\,|\widehat f(l)|^2\,\binom{l+2\lambda-1}{l}\\
&=2\,(2\lambda+1)\sum_{l=1}^{m}(l+\lambda)\,\binom{l+2\lambda}{l-1}\\
&=\frac{2\,(2m+2\lambda+1)\,(m+2\lambda+1)!\,(2\lambda+1)\,(\lambda+1)}{(m-1)!\,(2\lambda+3)!}
\end{align*}
compare Lemma~\ref{lem:Uf2}. Combining these expressions together according to~\eqref{eq:varS} and~\eqref{eq:varM} yields~\eqref{eq:varSPhim} and~\eqref{eq:varMPhim}, respectively.
\hfill$\Box$

Values of space variance, momentum variance and uncertainty product for several values of~$\lambda$ and~$m$ are collected in the following table. Note that the lower bound for~$U(f)$ is equal to~$\lambda+\frac{1}{2}$.\\\renewcommand{\baselinestretch}{1.6}
\begin{center}\begin{small}
\begin{tabular}{||c||c|c|c|c|c|c||}
\hline\hline
$\begin{array}{c}\text{var}_S(\Phi_m)\\[-0.8em]\text{var}_M(\Phi_m)\\[-0.8em]U(\Phi_m)\end{array}$&$\lambda=\frac{1}{2}$&$\lambda=1$
   &$\lambda=\frac{3}{2}$&$\lambda=2$&$\lambda=\frac{5}{2}$&$\lambda=3$\\\hline\hline
m=1&$\begin{array}{c}3\\[-0.8em]1.5\\[-0.8em]2.12\end{array}$&$\begin{array}{c}5.25\\[-0.8em]2.4\\[-0.8em]3.55\end{array}$
&$\begin{array}{c}8\\[-0.8em]3.33\\[-0.8em]5.16\end{array}$&$\begin{array}{c}11.3\\[-0.8em]4.29\\[-0.8em]6.94\end{array}$
&$\begin{array}{c}15\\[-0.8em]5.25\\[-0.8em]8.87\end{array}$&$\begin{array}{c}19.3\\[-0.8em]6.22\\[-0.8em]10.94\end{array}$
\\\hline
m=2&$\begin{array}{c}1.25\\[-0.8em]4\\[-0.8em]2.24\end{array}$&$\begin{array}{c}2.06\\[-0.8em]6\\[-0.8em]3.52\end{array}$
&$\begin{array}{c}3\\[-0.8em]8\\[-0.8em]4.9\end{array}$&$\begin{array}{c}4.06\\[-0.8em]10\\[-0.8em]6.37\end{array}$
&$\begin{array}{c}5.25\\[-0.8em]12\\[-0.8em]7.94\end{array}$&$\begin{array}{c}6.56\\[-0.8em]14\\[-0.8em]9.59\end{array}$
\\\hline
m=3&$\begin{array}{c}0.78\\[-0.8em]7.5\\[-0.8em]2.42\end{array}$&$\begin{array}{c}1.25\\[-0.8em]10.8\\[-0.8em]3.67\end{array}$
&$\begin{array}{c}1.78\\[-0.8em]14\\[-0.8em]4.99\end{array}$&$\begin{array}{c}2.36\\[-0.8em]17.14\\[-0.8em]6.36\end{array}$
&$\begin{array}{c}3\\[-0.8em]20.25\\[-0.8em]7.79\end{array}$&$\begin{array}{c}3.69\\[-0.8em]23.33\\[-0.8em]9.29\end{array}$
\\\hline
m=4&$\begin{array}{c}0.56\\[-0.8em]12\\[-0.8em]2.6\end{array}$&$\begin{array}{c}0.89\\[-0.8em]16.8\\[-0.8em]3.87\end{array}$
&$\begin{array}{c}1.25\\[-0.8em]21.33\\[-0.8em]5.16\end{array}$&$\begin{array}{c}1.64\\[-0.8em]25.71\\[-0.8em]6.5\end{array}$
&$\begin{array}{c}2.06\\[-0.8em]30\\[-0.8em]7.87\end{array}$&$\begin{array}{c}2.52\\[-0.8em]34.22\\[-0.8em]9.28\end{array}$
\\\hline
m=5&$\begin{array}{c}0.44\\[-0.8em]17.5\\[-0.8em]2.78\end{array}$&$\begin{array}{c}0.69\\[-0.8em]24\\[-0.8em]4.07\end{array}$
&$\begin{array}{c}0.96\\[-0.8em]30\\[-0.8em]5.37\end{array}$&$\begin{array}{c}1.25\\[-0.8em]35.71\\[-0.8em]6.68\end{array}$
&$\begin{array}{c}1.56\\[-0.8em]41.25\\[-0.8em]8.02\end{array}$&$\begin{array}{c}1.89\\[-0.8em]46.67\\[-0.8em]9.39\end{array}$
\\\hline
m=6&$\begin{array}{c}0.36\\[-0.8em]24\\[-0.8em]2.94\end{array}$&$\begin{array}{c}0.56\\[-0.8em]32.4\\[-0.8em]4.27\end{array}$
&$\begin{array}{c}0.78\\[-0.8em]40\\[-0.8em]5.58\end{array}$&$\begin{array}{c}1.01\\[-0.8em]47.14\\[-0.8em]6.89\end{array}$
&$\begin{array}{c}1.25\\[-0.8em]54\\[-0.8em]8.22\end{array}$&$\begin{array}{c}1.51\\[-0.8em]60.67\\[-0.8em]9.56\end{array}$
\\\hline
m=7&$\begin{array}{c}0.31\\[-0.8em]31.5\\[-0.8em]3.11\end{array}$&$\begin{array}{c}0.47\\[-0.8em]42\\[-0.8em]4.46\end{array}$
&$\begin{array}{c}0.65\\[-0.8em]51.33\\[-0.8em]5.79\end{array}$&$\begin{array}{c}0.84\\[-0.8em]60\\[-0.8em]7.11\end{array}$
&$\begin{array}{c}1.04\\[-0.8em]68.25\\[-0.8em]8.43\end{array}$&$\begin{array}{c}1.25\\[-0.8em]76.22\\[-0.8em]9.76\end{array}$
\\\hline
m=15&$\begin{array}{c}0.14\\[-0.8em]127.5\\[-0.8em]4.19\end{array}$&$\begin{array}{c}0.21\\[-0.8em]162\\[-0.8em]5.83\end{array}$
&$\begin{array}{c}0.28\\[-0.8em]190\\[-0.8em]7.35\end{array}$&$\begin{array}{c}0.36\\[-0.8em]214.29\\[-0.8em]8.8\end{array}$
&$\begin{array}{c}0.44\\[-0.8em]236.25\\[-0.8em]10.2\end{array}$&$\begin{array}{c}0.52\\[-0.8em]256.67\\[-0.8em]11.57\end{array}$
\\\hline
m=31&$\begin{array}{c}0.07\\[-0.8em]511.5\\[-0.8em]5.79\end{array}$&$\begin{array}{c}0.1\\[-0.8em]632.4\\[-0.8em]7.92\end{array}$
&$\begin{array}{c}0.13\\[-0.8em]723.33\\[-0.8em]9.82\end{array}$&$\begin{array}{c}0.17\\[-0.8em]797.14\\[-0.8em]11.57\end{array}$
&$\begin{array}{c}0.2\\[-0.8em]860.25\\[-0.8em]13.21\end{array}$&$\begin{array}{c}0.24\\[-0.8em]916.22\\[-0.8em]14.78\end{array}$
\\\hline
m=63&$\begin{array}{c}0.03\\[-0.8em]2047.5\\[-0.8em]8.01\end{array}$&$\begin{array}{c}0.05\\[-0.8em]2494.8\\[-0.8em]10.99\end{array}$
&$\begin{array}{c}0.06\\[-0.8em]2814\\[-0.8em]13.47\end{array}$&$\begin{array}{c}0.08\\[-0.8em]3060\\[-0.8em]15.74\end{array}$
&$\begin{array}{c}0.1\\[-0.8em]3260.25\\[-0.8em]17.83\end{array}$&$\begin{array}{c}0.11\\[-0.8em]3430\\[-0.8em]19.79\end{array}$
\\\hline
m=127&$\begin{array}{c}0.02\\[-0.8em]8191.5\\[-0.8em]11.38\end{array}$&$\begin{array}{c}0.02\\[-0.8em]9906\\[-0.8em]15.34\end{array}$
&$\begin{array}{c}0.03\\[-0.8em]11091.3\\[-0.8em]18.76\end{array}$&$\begin{array}{c}0.04\\[-0.8em]11974.3\\[-0.8em]21.82\end{array}$
&$\begin{array}{c}0.05\\[-0.8em]12668.3\\[-0.8em]24.61\end{array}$&$\begin{array}{c}0.06\\[-0.8em]13236.2\\[-0.8em]27.2\end{array}$
\\\hline
m=255&$\begin{array}{c}0.01\\[-0.8em]32767.5\\[-0.8em]16.05\end{array}$&$\begin{array}{c}0.01\\[-0.8em]39474\\[-0.8em]21.58\end{array}$
&$\begin{array}{c}0.02\\[-0.8em]44030\\[-0.8em]26.33\end{array}$&$\begin{array}{c}0.02\\[-0.8em]47357.1\\[-0.8em]30.55\end{array}$
&$\begin{array}{c}0.02\\[-0.8em]49916.3\\[-0.8em]34.37\end{array}$&$\begin{array}{c}0.03\\[-0.8em]51963\\[-0.8em]37.9\end{array}$
\\\hline\hline
\end{tabular}\\
\end{small}\end{center}
\renewcommand{\baselinestretch}{1}\pagebreak
\begin{lem}The space variance, momentum variance and uncertainty product of~$\Phi_m$ behave for~$m\to\infty$ like
\begin{align*}
\text{var}_S(\Phi_m)&\sim\frac{2\lambda+1}{m},\\
\text{var}_M(\Phi_m)&\sim\frac{2\lambda+1}{2\lambda+3}\cdot m^2,\\
U(\Phi_m)&\sim\frac{2\lambda+1}{\sqrt{2\lambda+3}}\cdot m^{1/2}.
\end{align*}
\end{lem}

\begin{bfseries}Proof. \end{bfseries}By elementary calculations.\hfill$\Box$

\begin{bfseries}Remark. \end{bfseries}For~$\lambda=\frac{1}{2}$ this coincides with the result obtained by Fern\'andez in~\cite[Section~5.2]{nFPhD}.

In a similar way, we can describe behavior of~$\Phi_m$'s for growing space dimension, as it is described in the next lemma.

\begin{lem}The space variance, momentum variance and uncertainty product of~$\Phi_m$ behave for~$\lambda\to\infty$ like
\begin{align*}
\text{var}_S(\Phi_m)&\sim\begin{cases}4\lambda^2&\text{for }m=1,\\\frac{1}{m^2}\cdot\lambda^2&\text{for }m>1,\end{cases}\\
\text{var}_M(\Phi_m)&\sim2m\lambda,\\
U(\Phi_m)&\sim\begin{cases}(2\lambda)^{3/2}&\text{for }m=1,\\\left(\frac{2}{m}\right)^{1/2}\cdot(\lambda)^{3/2}&\text{for }m>1.\end{cases}
\end{align*}
\end{lem}

\begin{bfseries}Proof. \end{bfseries}By elementary calculations.\hfill$\Box$

\section{Spherical basis functions}\label{sec:SBA}

\begin{df}
A continuous function $G:\,[-1,1]\to\mathbb R$ is called positive definite on~$\mathcal S^n$ if for every $L\in\mathbb N$ and any sequence of points $(x_l)_{l=1}^L\subseteq\mathcal S^n$ the corresponding Gramian matrix $A=\left(G(x_j\cdot x_k)\right)_{\genfrac{}{}{0pt}{}{j=1,\dots,L}{k=1,\dots,L}}$ is positive semidefinite. Further, if~$A$ is positive definite, then~$G$ is called strictly positive definite.
\end{df}

Positive definite functions and strictly positive definite functions are used for interpolation scattered data on a sphere, and strict positive definiteness is much more advantageous  for stability of algorithms, compare~\cite{XC92}, \cite{RS96}, and~\cite{CMS03} for the discussion. The starting point in determining whether a function is positive definite or not is the classical result obtained by Schoenberg in~\cite{iS42}.

\begin{thm}\label{thm:schoenberg}A continuous function $G:\,[-1,1]\to\mathbb R$ is positive semidefinite on~$\mathcal S^n$ if and only if it has the form
$$
G=\sum_{l=0}^\infty a_lC_l^\lambda
$$
with $a_l\geq0$ and $\sum_{l=0}^\infty C_l^\lambda(1)<\infty$.
\end{thm}

Further, for positive definiteness of cardinality~$L$ (i.e., for some fixed~$L$ instead of $L\in\mathbb N$, compare~\cite{RS96}) we have the following sufficient condition, proven in~\cite{XC92}.

\begin{thm}Let~$G$ satisfy conditions of Theorem~\ref{thm:schoenberg} and let $(x_l)_{l=1}^L\subseteq\mathcal S^n$. In order that the matrix $A=\left(G(x_j\cdot x_k)\right)_{\genfrac{}{}{0pt}{}{j=1,\dots,L}{k=1,\dots,L}}$ is positive definite it is sufficient that the coefficients~$a_l$ be positive for $0\leq l<L$.
\end{thm}

Finally, a characterization of strict positive definiteness is given in~\cite{CMS03}.

\begin{thm}Let~$G$ satisfy conditions of Theorem~\ref{thm:schoenberg}. In order that~$G$ is strictly positive definite on~$\mathcal S^n$ it is necessary and sufficient that infinitely many coefficients~$a_l$ with odd index~$l$ and infinitely many coefficients~$a_l$ with even index~$l$ be positive.
\end{thm}

According to these theorems, weighted scaling functions~$\phi_j$ are positive definite, strictly positive definite of cardinality~$2^{j-1}$, but they are not strictly positive definite.

\section{Appendix}
\begin{lem}\label{lem:dim_Pi}For $n\geq2$ it holds
$$
\sum_{l=0}^m\frac{(2l+n-1)(l+n-2)!}{(n-1)!\,l!}=\frac{(n+2m)(n+m-1)!}{n!\,m!}.
$$
\end{lem}

\begin{bfseries}Proof. \end{bfseries}(by induction) It is straightforward to verify that the formula holds for $m=0$. Suppose, it is valid for some~$m$. For $m+1$ we have
\begin{align*}
\sum_{l=0}^{m+1}&\frac{(2l+n-1)(l+n-2)!}{(n-1)!\,l!}\\&=\frac{(n+2m)(n+m-1)!}{n!\,m!}+\frac{(n+2m+1)(n+m-1)!}{(n-1)!\,(m+1)!}\\
&=\left[(m+1)(n+2m)+n(n+2m+1)\right]\cdot\frac{(n+m-1)!}{n!\,(m+1)!}\\
&=(n+2m+2)(n+m)\cdot\frac{(n+m-1)!}{n!\,(m+1)!}=\frac{(n+2m+2)(n+m)!}{n!\,(m+1)!}.
\end{align*}
\hfill$\Box$

\begin{lem}\label{lem:Uf2}
For $\lambda>0$ we have
$$
\sum_{l=1}^m(l+\lambda)\,\binom{l+2\lambda}{l-1}=\frac{(2m+2\lambda+1)\,(m+2\lambda+1)!\,(\lambda+1)}{(m-1)!\,(2\lambda+3)!}.
$$
\end{lem}

\begin{bfseries}Proof. \end{bfseries}(by induction) It is straightforward to verify that the formula holds for $m=1$. Suppose, it is valid for some~$m$. For $m+1$ we have
\begin{align*}
\sum_{l=1}^{m+1}&\,(l+\lambda)\,\binom{l+2\lambda}{l-1}\\
&=\frac{(2m+2\lambda+1)\,(m+2\lambda+1)!\,(\lambda+1)}{(m-1)!\,(2\lambda+3)!}+\frac{(m+\lambda+1)\,(m+2\lambda+1)!}{m!\,(2\lambda+1)!}\\
&=\frac{(m+2\lambda+1)!\left[m\,(2m+2\lambda+1)\,(\lambda+1)+(m+\lambda+1)\,(2\lambda+2)\,(2\lambda+3)\right]}{m!\,(2\lambda+3)!}\\
&=\frac{(m+2\lambda+2)!\,(2m+2\lambda+3)\,(\lambda+1)}{m!\,(2\lambda+3)!}.
\end{align*}
\hfill$\Box$

\end{document}